\documentclass{amsart}
\usepackage{graphicx, amssymb, mathrsfs}
\usepackage[all, cmtip]{xy}
\newtheorem{thm}{Theorem}[section]

\newtheorem{lem}[thm]{Lemma}

\theoremstyle{definition}

\theoremstyle{remark}

\numberwithin{equation}{section}
\newcommand{\norm}[1]{\left\Vert#1\right\Vert}
\newcommand{\abs}[1]{\left\vert#1\right\vert}
\newcommand{\set}[1]{\left\{#1\right\}}

\newcommand{\To}{\longrightarrow}

\newcommand{\cplx}{\mathbb{C}}
\newcommand{\proj}{\mathbb{P}}

\newcommand{\Aut}{\rm{Aut}}

\begin{document}

\title[Cut-off Function Lemma in $\proj^k$]{Cut-off Function Lemma in $\proj^k$}%
\author{Taeyong Ahn}%
\address{5840 East Hall, 530 church street, Ann Arbor, MI 48109-1043}%
\email{taeyonga@umich.edu}%

\thanks{TBA}%
\subjclass{TBA}%

\begin{abstract}
In this note, we compute a cut-off function over $\proj^k$. Let sufficiently small $\delta>0$ be given. When we are given a compact set $K$ in $\proj^k$ and a prescribed open neighborhood $K_\delta$ of $K$, we find a smooth cut-off function $\chi_\delta$ such that $\chi_\delta\equiv 1$ over $K$ and $supp(\chi_\delta)\subseteq K_\delta$, where $K_\delta$ denotes the set of points whose distance to $K$ is less than $\delta$ with respect to the Fubini-Study metric of $\proj^k$. Moreover, we estimate the bound of the derivatives of $\chi_\delta$ in terms of $\delta$. It seems to be well-known, but we want to provide detailed computations. They are very elementary.
\end{abstract}

\maketitle
\section{Introduction}

In this note, our space is $\proj^k$ and we assume that the distance is measured with respect to the Fubini-Study metric if we do not specify.\\

Let $\delta_0>0$ be given. We consider $0<\delta<\delta_0$. Let $K\subseteq\proj^k$ be compact and $K_\delta$ a $\delta$-neighborhood of $K$, that is, the set of points whose distance to $K$ is less than $\delta$ with respect to the Fubini-Study metric. We want to prove the following lemma:
\begin{lem}\label{lem:main}
There exists a smooth cut-off function $\chi_\delta:\proj^k\to [0, 1]$ such that $\chi_\delta\equiv 1$ over $K$ and $supp(\chi_\delta)\subseteq K_\delta$. Moreover, $\norm{\chi_\delta}_{C^\alpha}\lesssim \abs{\delta}^{-\alpha}$ as $\delta$ varies.
\end{lem}
Here, $\norm{\cdot}_{C^\alpha}$ denotes the $C^\alpha$-norm of the function. The idea is simply to smooth out a characteristic function by convolution (of the Lie group of automorphisms over $\proj^k$).

\section{Family of Local Coordinate Charts of $\proj^k$}
It suffices to prove the lemma for a fixed family of local coordinate charts. Thus, we will fix one as follows.\\

For $\proj^k$, we can find $k$ natural affine coordinate charts covering $\proj^k$ of the form $\set{[z_0: ...: z_{i-1}: 1: z_{i+1}: ...: z_k]| z_j\in\cplx \textrm{ for }j\neq i}$ for $i=0, ..., k$, which we will call the $Z_i$-coordinate chart. For this chart, there is a natural coordinate map $\zeta_i:Z_i\to\cplx^{i}\times\set{1}\times\cplx^{k-i}$ defined by $\zeta_i([z_0: ...: z_{i-1}: 1: z_{i+1}: ...: z_k])=(z_0, ..., z_{i-1}, 1, z_{i+1}, ..., z_k)$.\\

We defined a norm $\norm{\cdot}_i$ defined by
\begin{displaymath}
\norm{(z_0, ..., z_{i-1}, 1, z_{i+1}, ..., z_k)}_i=(\abs{z_0}^2+ ...+ \abs{z_{i-1}}^2+ \abs{z_{i+1}}^2+ ...+ \abs{z_k}^2)^\frac{1}{2}
\end{displaymath}
for each $\cplx^{i}\times\set{1}\times\cplx^{k-i}$.

\section{Automorphism group of $\proj^k$}\label{sec:Auto}
The group $\Aut(\proj^k)=PGL(k+1, \cplx)$ of automorphisms of $\proj^k$ is a complex Lie group of complex dimension $k^2+2k$. An element of $\Aut(\proj^k)$ can be understood as an equivalence class of the complex $(k+1)\times(k+1)$ matrix group under the equivalence relation given by scaling.\\

Without loss of generality, we may consider a point $z\in Z_0$ and its coordinates $\zeta\in\set{1}\times\cplx^k$. Let $h=(0, h_1, h_2, ..., h_k)$ with $\abs{h_i}<\epsilon$ for sufficiently small $\epsilon>0$. Then $\zeta+h\in \set{1}\times\cplx^k$ is a very close point near $\zeta\in \set{1}\times\cplx^k$, where the addition is coordinatewise and we can find a unique linear map $G_h:\set{1}\times\cplx^k\to\set{1}\times\cplx^k$ defined by 
\begin{displaymath}
G_h=\left(\begin{array}{ccccc}
1 & 0 & 0 & \cdots & 0\\
h_1 & 1 & 0 & \cdots & 0\\
h_2 & 0 & 1 & \cdots & 0\\
\vdots & \vdots & \vdots & \ddots & \vdots\\
h_n & 0 & 0 & \cdots & 1
\end{array}\right)
\end{displaymath}
 such that $G_h(\zeta)=\zeta+h$. Note that $G_h\circ G_{-h}=G_{-h}\circ G_h=Id$.\\

Using the exponential map of Lie algebra to Lie group, we can find holomorphic coordinates $\psi:sl(k+1, \cplx)\to PGL(k+1, \cplx)$ near $Id\in PGL(k+1, \cplx)$ where $sl(k+1, \cplx)$ is the special linear Lie algebra, which is the set of $(k+1)\times(k+1)$ matrices with zero trace. Near the $Id\in PGL(k+1, \cplx)$, we can also find a representation $PGL(k+1, \cplx)\to GL(k+1, \cplx)$ by picking a $(k+1)\times (k+1)$ matrix with the $(1, 1)$-component being $1$. Let $\phi$ denote this representation. We consider the following diagram

\begin{displaymath}
\begin{array}{ccc}
sl(k+1, \cplx) & \To^{H_h} & sl(k+1, \cplx)\\
\psi\downarrow & & \psi\downarrow\\
PGL(k+1, \cplx) & \To^{[G_h]} & PGL(k+1, \cplx)\\
\phi\downarrow & & \phi\downarrow\\
GL(k+1, \cplx) & \To^{\overline{G_h}} & GL(k+1, \cplx),
\end{array}
\end{displaymath}
where in the second line, $[\cdot]$ means the equivalence class that contains the inside element,  $[G_h][A]=[A\cdot G_h]$ for $[A]\in PGL(k+1, \cplx)$, and $\overline{G_h}$ is defined as follows:
\begin{eqnarray*}
&\left(\begin{array}{cccc}
1 & a_{1, 2} & \cdots & a_{1, k+1}\\
a_{2, 1} & a_{2, 2} & \cdots & a_{2, k+1}\\
a_{3, 1} & a_{3, 2} & \cdots & a_{3, k+1}\\
\vdots & \vdots & \ddots & \vdots\\
a_{k+1, 1} & a_{k+1, 2} & \cdots & a_{k+1, k+1}
\end{array}\right)&\\
&\downarrow \overline{G_h}&\\
&\frac{1}{1+\sum_{i=2}^{k+1}a_{2, i}\cdot h_{i-1}}\left(\begin{array}{cccc}
1 & a_{1, 2} & \cdots & a_{1, k+1}\\
a_{2, 1}+\sum_{i=2}^{k+1}a_{2, i}\cdot h_{i-1} & a_{2, 2} & \cdots & a_{2, k+1}\\
a_{3, 1}+\sum_{i=2}^{k+1}a_{3, i}\cdot h_{i-1} & a_{3, 2} & \cdots & a_{3, k+1}\\
\vdots & \vdots & \ddots & \vdots\\
a_{k+1, 1}+\sum_{i=2}^{k+1}a_{k+1, i}\cdot h_{i-1} & a_{k+1, 2} & \cdots & a_{k+1, k+1}
\end{array}\right)&.
\end{eqnarray*}
Note that $H_h, [G_h]$ and $\overline{G_h}$ in the diagram are not defined over the entire space. However, there exists a sufficiently small $\epsilon>0$ such that for all $\set{h_i}_{i=1}^n$ with $\abs{h_i}<\epsilon$ for $i=1, ..., n$, $\overline{G_h}$ is well-defined over all $A\in GL(n+1, \cplx)$ with $\norm{A-Id}<\epsilon$ and with the $(1, 1)$-component of $A$ being $1$, where $\norm{\cdot}$ is the standard matrix norm. Since $\phi$ and $\psi$ are local biholomorphisms, we can also find correspoinding subsets in $PGL(k+1, \cplx)$ and $sl(k+1, \cplx)$.\\

We identify $sl(k+1, \cplx)$ with $\cplx^{k^2+2k}$ and the set of representations of $PGL(k+1, \cplx)$ with $\cplx^{k^2+2k}$. For convenience, we use $x=(x_1, ..., x_{k^2+2k})$ for $sl(k+1, \cplx)$ and $\xi=(\xi_1, ..., \xi_{k^2+2k})$ for the other. Then
\begin{displaymath}
\begin{array}{ccc}
\cplx^{k^2+2k} & \To^{H_h} & \cplx^{k^2+2k}\\
\psi\downarrow & & \psi\downarrow\\
PGL(k+1, \cplx) & \To^{[G_h]} & PGL(k+1, \cplx)\\
\phi\downarrow & & \phi\downarrow\\
\cplx^{k^2+2k} & \To^{\overline{G_h}} & \cplx^{k^2+2k}.
\end{array}
\end{displaymath}

We denote $\phi\circ\psi$ by $\Phi$. Then, $\xi_i=\Phi_i(x_1, ..., x_{k^2+2k})$ for $i=1, ..., k^2+2k$ and the map $H_h=\Phi^{-1}\circ \overline{G_h}\circ \Phi$ is a map from $\cplx^{k^2+2k}$ to $\cplx^{k^2+2k}$. Note that in our case, $\psi, \phi$ are smooth and $\overline{G_h}$ is smooth with respect to $h$.

\section{Measures on $sl(k+1, \cplx)$}
Recall that $x$ is used for $sl(k+1, \cplx)$. Let $\lambda$ denote the standard Euclidean measure on $sl(k+1, \cplx)$. We assign the standard matrix norm $\norm{x}_s$ to each $x\in sl(k+1, \cplx)$. We consider a smooth radial probability measure $\mu$ over the coordinate $sl(k+1, \cplx)$ centered at $O\in sl(k+1, \cplx)$ with its support $\norm{x}_s<\sigma$ for sufficiently small $\sigma>0$, which makes $\Phi(\set{\norm{x}_s<\sigma})\subseteq \set{\norm{A-Id}<\epsilon}$. Then, $d\mu=M(x)d\lambda$ where $M$ is a smooth function defined on $sl(k+1, \cplx)$ and has support in $\norm{x}_s<\sigma$.\\

Let $h_\theta: sl(k+1, \cplx)\To sl(k+1, \cplx)$ be a scaling map by $\theta$ for $\abs{\theta}\leq 1$. We define $\mu_\theta:=(h_\theta)_*(\mu)$. Then, $\mu_\theta$ is a smooth measure for $\theta\neq 0$ and a Dirac measure at $O\in sl(k+1, \cplx)$ for $\theta=0$. Note that the support of $\mu$ is in $\set{\norm{x}_s\leq\theta\sigma}\subseteq \set{\norm{x}_s\leq\sigma}$.\\

For the better terminology, by the derivatives of $\mu_\theta$, we mean the derivatives of the Radon-Nikodym derivative of $\mu_\theta$ with respect to the standard Euclidean measure $\lambda$.

\section{Regularization}
In this section, we define a regularization of a bounded function and provide the estimate of the regularity.\\

Let $f$ be a bounded complex-valued function over $\proj^k$ with compact support. Without loss of generality, we may assume that $0\leq\abs{f}\leq 1$. Then, we define the $\theta$-regularization $f_\theta$ of $f$ as being
\begin{displaymath}
	f_\theta(z)=\int_{\Aut(\proj^k)}((\tau_x)_*f)(z)d\mu_\theta(x).
\end{displaymath}

Without loss of generality, we may assume that $z\in Z_0$. Let $\zeta\in\set{1}\times\cplx^k$ be the coordinates of $z$ and $F$ the representation of $f$ with respect to $\set{1}\times\cplx^k$. With respect to the coordinate $\set{1}\times\cplx^k$, we have the following representation:
\begin{eqnarray*}
F_\theta(\zeta+h)&=&\int_{sl(k+1, \cplx)}((\Phi(x))_*F)(G_h(\zeta))d\mu_\theta(x)\\
&=&\int_{sl(k+1, \cplx)}((\Phi(H_h(x)))_*F)(\zeta)d\mu_\theta(x)
\end{eqnarray*}
Note that $H_h$ is holomorphic and injective over the support of the measure $\mu_\theta$. By change of coordinates, we have
\begin{eqnarray*}
F_\theta(\zeta+h)&=&\int_{sl(k+1, \cplx)}((\Phi(H_h(x)))_*F)(\zeta)d\mu_\theta(x)\\
&=&\int_{sl(k+1, \cplx)}((\Phi(x))_*F)(\zeta)((H_h)_*d\mu_\theta)(x).
\end{eqnarray*}

With $\zeta$ fixed, the differentiation of the right hand side with respect to $h_i$'s makes sense since the measure is smooth. By the direct application of the definition of the derivative, the partial derivative of $F_\theta(\zeta)$ with respect to $\zeta_i$ at $\zeta$ is the same as the partial derivative of $F_\theta(\zeta+h)$ with respect to $h_i$ at $0$. Thus, we can see that $F_\theta$ is smooth. Moreover, we can estimate its regularity.\\

The $C^\alpha$-norm of $F_\theta$ completely depends on the value of $F$ near $\zeta$ and the derivatives of the measure with respect to $h$. It is not hard to see that $(H_h)_*[(h_\theta)_*d\lambda]=\abs{\theta}^{-2k^2-4k}d\lambda$. Indeed, $\Phi$ is a coordinate change map and $G_h$ is a linear shear map. Thus, it remains to estimate the $C^\alpha$-norm of $M$. So, since $(H_h)_*[(h_\theta)_*M]=M(\frac{1}{\theta}[\Phi^{-1}\circ\overline{G_h}\circ\Phi])$, the $C^\alpha$-norm of $(H_h)_*[(h_\theta)_*M]$ is bounded by the product of $\abs{\theta}^{-\alpha}$ and a constant multiple of $C^\alpha$-norms of $M$, $\Phi$ and $\Phi^{-1}$. Note that the latter is independent of $\theta$.\\

Putting all together, since $F$ is bounded, the support of the measure is $\norm{x}\leq \theta\sigma$ and $dim_{\cplx}sl(k+1, \cplx)$ is $k^2+2k$,
\begin{equation}\label{eqn:estimate}
{f_\theta}_{C^{\alpha}}\lesssim\abs{\theta}^{-2k^2-4k-\alpha}\abs{\theta}^{2k^2+4k}\norm{f}_{C^{\alpha}}=\abs{\theta}^{-\alpha}\norm{f}_{C^{\alpha}}.
\end{equation}
Note that it can be more precise when we estimate the absolute value at a point in terms of its neighborhood with compact closure.



\section{Main Cut-off Function Lemma}
We consider two kinds of open balls in $\set{1}\times\cplx^k$. One is induced from the Fubini-Study metric of $\proj^k$ and the other is from the standard Euclidean metric $\norm{\cdot}_0$. The open ball centered at $\zeta\in\set{1}\times\cplx^k$ and of radius $r>0$ of first kind is denoted by $B_F(\zeta, r)$ and that of second kind is denoted by $B_E(\zeta, r)$. Then, by comparison of the infinitesimal versions of the two metrics, we know that $B_E(\zeta, \frac{r}{2}\norm{\zeta}_0)\subseteq B_F(\zeta, r)$.\\

\begin{proof}[The proof of Lemma ~\ref{lem:main}]
Note that $\Phi$ is holomorphic near the closure of the neighborhood of $\set{\norm{x}_s<\sigma}$, we can find a constant $C>0$ such that $\frac{1}{C}\norm{\Phi(x)-Id}<\norm{x}_s<C\norm{\Phi(x)-Id}$ for $\set{\norm{x}_s<\sigma}$. Here, $C$ is independent of $\delta$ and $\theta$. Recall that $\norm{\Phi(x)(\zeta)-\zeta}_0\leq\norm{\Phi(x)-Id}\norm{\zeta}_0$. We take a $\theta$ such that $\abs{\theta}\leq 1$ and such that $C\theta\sigma\leq\frac{\delta_0}{4}$. Let $C':=\frac{C\theta\sigma}{\delta_0/4}\leq 1$. Then, for all $0<\delta<\delta_0$, we take its corresponding $\theta$ to satisfy $C\theta\sigma=C'\frac{\delta}{4}$. Note that This $C'$ is fixed with respect to $\theta$ and $\delta$. Then, for each $0<\delta<\delta_0$ and for its $\theta$, we have that for $\set{\norm{x}_s<\sigma}$,
\begin{eqnarray}\label{eqn:domain}
\norm{\Phi(x)(\zeta)-\zeta}_0&\leq&\norm{\Phi(x)-Id}\norm{\zeta}_0\leq C\norm{x}_s\norm{\zeta}_0\\
\nonumber &\leq& C\theta\sigma\norm{\zeta}_0 =\frac{C'\delta}{2}\frac{\norm{\zeta}_0}{2}\leq \frac{\delta}{2}\frac{\norm{\zeta}_0}{2}.
\end{eqnarray}
\\

Consider $K\subseteq K_{\frac{\delta}{2}}\subseteq K_\delta$. Let $\chi_K$ be the characteristic function whose support is exactly $K_{\frac{\delta}{2}}$. Then $(\chi_K)_\theta$ is the desired function with the desired estimate. Indeed, the estimate is straight forward by plugging-in $C\theta\sigma=C'\frac{\delta}{4}$ into Estimate ~\ref{eqn:estimate}. Equation ~\ref{eqn:domain} proves the support of the function and its region over which the function is identically $1$.\\

So far, we have considered over $Z_0$ only. The above argument can be directly applied to each $Z_i$ for $i=0, ..., k$ in the exactly same way. Indeed, we use the same measure on $\Aut(\proj^k)$ and the same constants $C, C'$ and $\theta$ to $Z_i$ for $i=1, ..., k$ as in the case of $Z_0$. Thus, we have just proved the lemma.\\
\end{proof}
\end{document}